
\input amstex
\documentstyle{amsppt}
\magnification \magstep1
\parskip10pt
\parindent.3in
\pagewidth{5.7in}
\pageheight{7.2in}
\NoRunningHeads

\centerline {\bf An algebraic version of Demailly's asymptotic Morse
Inequalities.}
\vskip.4cm
\centerline {\smc Flavio Angelini}
\vskip.4cm
\topmatter
\abstract We give an elementary algebraic proof of some asymptotic
estimates (called by Demailly asymptotic Morse inequalities) for the
dimensions of cohomology groups of the difference
of two ample line bundles on a smooth complex projective variety of
any dimension.
\endabstract
\subjclass Primary 14F99, Secondary 32J99
\endsubjclass
\endtopmatter
\document

\def\qsum{\sum_{i=0}^q}
\def\q-1sum{\sum_{i=0}^{q-1}}
\def\jsum{\sum_{l=0}^{j-1}}

\def\0{\Cal O}
\def\K{\chi}

\vskip.4cm
\heading 1. Introduction
\endheading

\vskip.4cm
In \cite{De1} Demailly proves the asymptotic Morse Inequalities for a
holomorphic
line bundle on a compact K\"ahler Manifold. His methods are
essentially analytic and his first formulation involves some curvature
integral in general quite difficult to compute.
A formulation of these inequalities in a more algebraic setting appears
in \cite{De2, \S7}.
Motivated by the work of Siu on an effective version of Matsusaka's big theorem
\cite{Siu} and by a work of Trapani on the difference of two ample line bundles
\cite{Tra}, Demailly deduces, from his Morse inequalities, some asymptotic
estimates for the dimensions of cohomology groups
of the difference of two numerically effective line bundles in which only
algebraic invariants are involved \cite{De3, \S12}.
The purpose of this note is to give an elementary algebraic proof
of these estimates in the case of a smooth projective variety.

The theorem is the following:

\vskip.2cm
\proclaim {Theorem A}
Let $L=F-G$ be a line bundle on a smooth complex projective variety $X$ of
dimension
$n$, where $F$ and $G$ are ample line bundles. Then, for
every $q=0,1...,n$ and $k$ sufficiently large, we have:

$$\align
h^q(X,kL)&\leq k^n \frac{F^{n-q}\cdot G^q}{(n-q)!q!} + o(k^n) \tag1 \\
\qsum (-1)^{q-i}h^i(X,kL)&\leq \frac{k^n}{n!} \qsum (-1)^{q-i} \binom{n}{i}
     F^{n-i}\cdot G^i +o(k^n) \tag2 \\
\endalign \nopagebreak$$

Demailly \cite{De3, \S12} calls (1) and (2) respectively weak and strong
asymptotic Morse
inequalities.
\endproclaim

\vskip.2cm
The case $q=1$ gives a lower bound for the dimension of the space of
sections $h^0(X,kL)$ which was used by Siu \cite{Siu} in the course of his
proof of an effective version of
Matsusaka's big theorem.
An algebraic proof in this case was given by R. Lazarsfeld and F. Catanese
\cite{De3, \S12}.

More generally the theorem holds under milder positivity hypothesis on $F$,
depending on $q$, i.e. when $F$ is ample off a set of dimension less or equal
than $q-1$.

\proclaim{Theorem B}
Suppose $G$ is an ample line bundle and $F$ is a line bundle having the
property that there exist an ample vector
bundle $E$ and a surjective map
$$E \longrightarrow \Cal I_{V} \otimes F \longrightarrow0$$
where $\Cal I_{V}$ is the ideal sheaf of a subset $V$.
Then, for $k$ sufficiently large and $q\geq dim(V)+1$, the inequalities (1) and
(2) of
Theorem A hold for $L=F-G$.
\endproclaim

\vskip.2cm
Clearly Theorem B implies Theorem A.

\vskip.2cm
As a consequence we can obtain asymptotic estimates for the dimensions
of the cohomology groups of a big line bundle $L$.
If some multiple of $L$ is written as $aL=A+D$ where $A$ is an ample divisor
and $D$ an effective divisor, the upper bound is expressed in terms of the
first chern class of $L$,
the multiplicities of the singular points of $D$ and the curvature of the
tangent bundle $TX$ \cite{De2, \S7},\cite{De3, \S12}.

\vskip.4cm
This note is part of my Ph.D. thesis at UCLA, and I would like to
thank Rob Lazarsfeld for his guidance and encouragement.
I would like also to thank J.-P. Demailly for his prompt and detailed
correspondence about
the background on the topic.

\vskip.8cm
\heading 2. Proof of Theorem B
\endheading

\vskip.4cm
First we need the following:
\proclaim{Lemma}
Let $F$ be a line bundle on a smooth projective variety $X$ satisfying the
hypothesis of Theorem B. Then, for any very ample divisor $H$ on $X$, there
exists $k_o=k_o(X,F,H)$ such that
$$H^q(X,kF+aH)=0$$
for $q\geq \max(1,dim(V)+1)$, $k \geq k_o$ and any $a \geq 0$. \newline
In particular there exists $k_o=k_o(X,F)$ such that $H^q(X,kF)=0$ for
$k \geq k_o$ and \newline
$q \geq \max(1,dim(V)+1)$.
\endproclaim

\demo{Proof}
If $dim(V)=-1$, i.e. $V=$\O, then $F$ is ample. This is true also for
$dim(V)=0$
[For any ample $A$, $S^kE \otimes \0_X(-A)
$ is
globally generated for $k$ sufficiently large].
But then, in these two cases, there exists $k_o=k_o(X,F)$ such that
$k_oF=K_X+$(ample),
where $K_X$ is the canonical bundle of $X$.
Hence $kF+aH=K_X+$(ample) for $k \geq k_o$ and any $a \geq 0$ and the Lemma
follows
from Kodaira vanishing.

Now suppose $dim(V)=m$ with $m \geq 1$. Pick a general smooth irreducible
divisor $D$
in the linear series of $H$ that meets $V$ properly and consider, for $a \geq
1$,
 the short exact sequence
$$0\longrightarrow\0_X(kF+(a-1)H) \longrightarrow \0_X(kF+aH)\longrightarrow
 \0_D(kF+aH) \longrightarrow 0$$
By induction
$H^{q-1}(D,\0_D(kF+aH))=0$ for $k \geq k_o(D,F,H)$, $q \geq m+1$ and any $a
\geq 0$.
Hence $$H^q(X,\0_X(kF+(a-1)H))=H^q(X,\0_X(kF+aH))$$ for any $a \geq 0$ and $k
\geq k_o$.
But, for large $a$, $H^q(X,\0_X(kF+aH))=0$.  $\quad\quad\square$
\enddemo

To simplify notation let us set, for a divisor $B$
on $X$
$$\K_q(X,\0_X(B))= \qsum (-1)^{q-i}h^i(X,\0_X(B))$$
Observe that, if $$0 \rightarrow  A \rightarrow B \rightarrow C \rightarrow 0$$
  is an exact sequence of line bundles, then
$$\K_q(X,B) \leq \K_q(X,A)+ \K_q(X,C)$$
One can check that by writing the long exact sequence truncated
at the $q-th$ group. It follows that $\K_q$ satisfies similar subadditivity
with respect to filtrations.

\demo{Proof of Theorem B}
Let us pick a positive integer $a$ such that $aG$ is very ample and fix a
smooth
irreducible divisor $D$ in this linear series meeting $V$ properly.
In order to use induction on $dim(X)$ we will first prove an analogue of the
Theorem for the bundle $kF-jaG$ for some positive integer $j$. Subsequently
we'll deduce the stated inequalities for $L=F-G$.
Specifically, we'll prove:

$$\align h^q(X,(kF-jaG))&\leq k^{n-q} {(ja)}^{q} \frac{F^{n-q}\cdot
G^q}{(n-q)!q!} + o(k^{n-q})o(j^{q})\tag1 \\
\K_q(X,\0_X(kF-jaG))&\leq \frac{1}{n!} \qsum \left[ (-1)^{q-i}
k^{n-i}{(ja)}^{i} \binom{n}{i}
     F^{n-i} \cdot G^i +o(k^{n-i})o(j^{i}) \right]\tag2 \\
\endalign$$

We will only prove (2) being (1) similar and more simple.
To do so consider the exact sequence
$$0\longrightarrow\0_X(kF-jaG)\longrightarrow\0_X(kF)\longrightarrow
\0_{jD}(kF)\longrightarrow0\tag{$*$}$$
We have a filtration of $\0_{jD}$ whose quotients are line bundles of the form
$\0_D(-lD)=\0_D(-laG)$ with $l=0,\ldots,j-1$.
{}From the Lemma we have $H^q(X,\0_X(kF))=0$ for $k$ sufficiently large. Then
from ($*$) we get the long exact sequence in cohomology
$$0 \rightarrow H^o(X,\0_X(kF-jaG)) \rightarrow H^o(X,\0_X(kF)) \rightarrow
H^o(jD,\0_{jD}(kF)) \rightarrow \dots \nopagebreak$$
$$\rightarrow H^{q-1}(jD,\0_{jD}(kF)) \rightarrow
H^q(X,\0_X(kF-jaG)) \rightarrow H^q(X,\0_X(kF))=0$$
Therefore $$\K_q(X,\0_X(kF-jaG))=\K_q(X,\0_X(kF)) + \K_{q-1}(jD,\0_{jD}(kF))$$
Now from the filtration of $\0_{jD}$, we have
$$\K_{q-1}(jD,\0_{jD}(kF)) \leq \jsum \K_{q-1}(D,\0_D(kF-laG))$$
Since the bundle $\0_D(F)$ satisfies the hypothesis of the theorem with
$X$ replaced by $D$ and $V$ replaced by $V \cap D$, we can use induction
and assume:
$$\align &\K_{q-1}(D,\0_D(kF-laG)) \\ &\leq \frac{1}{(n-1)!} \q-1sum \left[
(-1)^{q-1-i} k^{n-1-i}l^i a^{i} \binom{n-1}i
F|_D^{n-1-i} \cdot G|_D^{i} + o(k^{n-1-i})o(l^{i}) \right] \endalign$$
Therefore, replacing $F|_D^{n-1-i} \cdot G|_D^{i}$ with $F^{n-1-i} \cdot
G^i \cdot D=aF^{n-1-i} \cdot G^{i+1}$, we find:
$$\align &\K_q(X,\0_X(kF-jaG))
 \leq \K_q(X,\0_X(kF))+ \\ &\frac{1}{(n-1)!}
\jsum \q-1sum \left[ (-1)^{q-1-i} k^{n-1-i}l^i a^{i+1} \binom{n-1}i
F^{n-1-i} \cdot G^{i+1} + o(k^{n-1-i})o(l^{i}) \right] \endalign$$
Using that $ \jsum l^i= \frac{j^{i+1}}{i+1}+o(j^{i+1})$ and
that, since $H^i(X,\0_X(kF))=0$ for $i \geq q$, \newline
$\K_q(X,\0_X(kF))=(-1)^q
\K(X,\0_X(kF))$ (where $\K$ is the usual characteristic of Euler-Poincar\'e)
the right hand side of the above inequality is equal to
$$\align &(-1)^q \left[ k^n \frac{F^n}{n!} + o(k^n) \right]+ \\
&\frac{1}{(n-1)!} \q-1sum \left[ (-1)^{q-1-i}
k^{n-1-i} \frac{j^{i+1}}{i+1} a^{i+1} \binom{n-1}i F^{n-1-i} \cdot
G^{i+1} + o(k^{n-1-i})o(j^{i+1}) \right] \endalign$$
Using the identity $ \frac{n}{i+1} \binom{n-1}i= \binom{n}{i+1}$ and
reindexing, the above can be written as
$$ \frac{1}{n!} \qsum \left[ (-1)^{q-i}k^{n-i}j^i a^i \binom ni F^{n-i} \cdot
G^i
+ o(k^{n-i})o(j^i) \right]$$
which is what we wanted.

\noindent
The first step of induction reduces, by cutting with hyperplane sections,
 to the case $q=0$ ($V=$\O) which works because
$$\K_0(X,\0_X(kF-laG))=h^0(X,\0_X(kF-laG)) \leq h^0(X,\0_X(kF))=k^n
\frac{F^n}{n!} +o(k^n)$$
To finish the proof of Theorem B let us take a large $k$ and write it as
$k=ja+a+r$ where $r<a$ and $a$ is chosen so that $mG$ is very ample
for any $m \geq a$.
Pick a smooth irreducible divisor $D \in |(a+r)G|$.
In this case we'll show (1), being (2) similar, but involving much more heavy
formulas. We have
$$\align &h^q(X,\0_X(k(F-G)) \leq h^q(X,\0_X(kF-jaG)) + h^{q-1}(D,\0_D(kF-jaG))
\\ &\leq k^{n-q}(ja)^q \frac {F^{n-q} \cdot G^q}{(n-q)!q!} +o(k^{n-q})o(j^q)
\\ &+k^{n-q}(ja)^{q-1}(a+r) \frac {F^{n-q} \cdot
G^q}{(n-q)!(q-1)!}+o(k^{n-q})o(j^{q-1}) \\
&\leq k^n \frac {F^{n-q} \cdot G^q}{(n-q)!q!}+o(k^n) \endalign$$
(2) works in the same way because
$$\K_q(X,\0_X(k(F-G)) \leq \K_q(X,\0_X(kF-jaG))+ \K_{q-1}(D,\0_D(kF-jaG))
\quad\quad\square$$
\enddemo

Now we would like to give an idea how one can use theorem B to
obtain similar results in the case of a big line bundle,
under some positivity hypothesis for the tangent bundle of
the variety. For more precise definitions and statements we
refer to \cite{De2, \S7} or \cite{De3, \S12}.

Let $X$ be a smooth projective variety, $L$ be a big line bundle
on $X$ and suppose there exists an ample line bundle $E$ such that
$TX \otimes \0_X(E)$ is ample.
 Some multiple of $L$ can be written as
$aL=A+D$ with $A$ ample and $D$ effective divisors and we can attatch
to $L$ some numbers $ 0=b_1 \leq \ldots \leq b_{n+1}$ obtained by counting
the multiplicities of the singular points of $D$.

If we let $J^r$ be the $r-th$ jet bundle associated to $aL$ then we get, for
$q=0, \ldots ,n$, a surjective map
$$J^{ab_{n-q+1}} \otimes \0_X(ab_{n-q+1}E) \longrightarrow
\Cal I_{V} \otimes \0_X(aL+ab_{n-q+1}E) \longrightarrow 0$$
where $V$ is a locus of singular points of $D$ that turns out,
by definition of the $b_i 's$, to be of dimension less or equal
than $q-1$.

Now $J^{ab_{n-q+1}} \otimes \0_X(ab_{n-q+1}E)$ is ample so we can apply
Theorem B to the bundles \newline $F=aL+ab_{n-q+1}E$ and $G=ab_{n-q+1}E$
and obtain asymptotic estimates for the dimensions of the cohomology
groups of $aL$. With further arguments we can actually obtain the same
results for $L$ itself.

\vskip4cm

\settabs\+ UCLA, Department of Mathematics \cr
\+ Department of Mathematics \cr
\+ UCLA \cr
\+ Los Angeles, CA 90024 \cr
\+ E-mail: fangelin\@math.ucla.edu \cr


\vskip2cm
\widestnumber\key{De3}

\Refs
\ref\key De1
    \by J.-P. Demailly
    \pages 189-229
    \paper Champs magn\'etiques et in\'egalit\'es de Morse pour la
d"-cohomologie
    \jour Ann. Inst. Fourier (Grenoble)
    \yr 1985 \vol 35
\endref

\ref\key De2
    \by J.-P. Demailly
    \paper Singular hermitian metrics on positive line bundles
    \paperinfo Conf. Complex algebraic varieties (Bayreuth, April 2-6, 1990)
    \inbook edited by K. Hulek, T. Peternell, M. Schneider, F. Schreyer
    \bookinfo Lecture Notes in Math., Vol. 1507
    \publ Springer-Verlag, Berlin, 1992.
\endref

\ref\key De3
    \by J.-P. Demailly
    \paper $L^2$ vanishing theorems for positive line bundles and adjunction
theory
    \paperinfo CIME session on Transcendental Methods in Alg. Geom.(Cetraro,
 Italy, July 1994)
    \jour Pr\'epublication de l'Inst. Fourier
    \yr 1994 \vol 288
\endref

\ref\key Siu
    \by Y. T. Siu
    \paper An effective Matsusaka big theorem
    \jour Ann. Inst. Fourier \vol 43 \yr 1993 \pages 1387-1405
\endref

\ref\key Tra
    \by S. Trapani
    \paper Numerical criteria for the positivity of the difference of
ample divisors
    \jour preprint
    \yr 1991
\endref
\endRefs
\enddocument